!Mode::"TeX:UTF-8"
\UseRawInputEncoding
\documentclass[review]{elsarticle}

\usepackage{lineno}
\usepackage[colorlinks,
            linkcolor=blue,       
            anchorcolor=blue,  
            citecolor=blue,        
            ]{hyperref}
\modulolinenumbers[5]

\usepackage{amssymb,amsmath}
\journal{}
\usepackage{geometry}
\usepackage{graphicx}
\usepackage{booktabs}
\usepackage{color}
\usepackage{colortbl}
\usepackage{epstopdf}
\usepackage{amssymb}

\newtheorem{definition}{Definition}
\newtheorem{theorem}{Theorem}
\newtheorem{remark}{Remark}

\usepackage{multirow}
\usepackage{graphicx}
\usepackage{subfigure}
\usepackage{cleveref}
\usepackage{titlesec}
\titlespacing*{\section}
{0pt}{0.5ex plus 1ex minus .1ex}{0.1ex plus .1ex}
\titlespacing*{\subsection}
{0pt}{0pt}{0.2pt}

\usepackage{bm}
\makeatletter
\newif\if@restonecol
\makeatother

\usepackage[linesnumbered,ruled,vlined]{algorithm2e}
\usepackage{algpseudocode}
\usepackage{amsmath}

\begin{document}
\begin{frontmatter}

\title{Generalized fractional grey system models: Memory effects perspective} 


\author[label1]{Wanli Xie}
 \address[label1]{Institute of EduInfo Science and Engineering, Nanjing Normal University, Nanjing 210097, China}

\author[label2]{Wen-Ze Wu\corref{cor1}}
\address[label2]{School of Economics and Business Administration, Central China Normal University, Wuhan 430079, China}
\cortext[cor1]{Corresponding author.}
\ead{wenzew@mails.ccnu.edu.cn}

\author[label3]{Chong Liu\corref{cor2}}
 \address[label3]{School of Science, Northeastern University, Shenyang 110819, China}
\cortext[cor2]{Corresponding author.}
\ead{liuchong@emails.imau.edu.cn}

\author[label4,label5]{Mark Goh}
\address[label4]{NUS Business School, National University of Singapore, 21 Lower Kent Ridge Road, Singapore}
\address[label5]{The Logistics Institute-Asia Pacific, National University of Singapore, 21 Lower Kent Ridge Road, Singapore}

\begin{abstract}
As an essential characteristics of fractional calculus, the memory effect is served as one of key factors to deal with diverse practical issues, thus has been received extensive attention since it was born. By combining the fractional derivative with memory effects and grey modeling theory, this paper aims to construct an unified framework for the commonly-used fractional grey models already in place. In particular, by taking different kernel and normalization functions, this framework can deduce some other new fractional grey models. To further improve the prediction performance, the four popular intelligent algorithms are employed to determine the emerging coefficients for the UFGM(1,1) model. Two published cases are then utilized to verify the validity of the UFGM(1,1) model and explore the effects of fractional accumulation order and initial value on the prediction accuracy, respectively. Finally, this model is also applied to dealing with two real examples so as to further explain its efficacy and equally show how to use the unified framework in practical applications.
\end{abstract}
\end{frontmatter}

\section{Introduction}
\label{sec:1}
Grey system theory, an effective approach to addressing the issues in uncertain systems, was proposed by Deng \cite{1}. The grey system theory can fully dig out the internal rule of systems with limited data, this enables it has very important application prospects and successfully solved many practical problems. Grey prediction model is the most essential model in the grey system theory, and it has become a research hotspot of scholars in recent years.
To be specific,
GM(1,1) \cite{2}, DGM(1,1) \cite{3} are two kinds of commonly-used grey models, there are many published papers focusing on them \cite{4,5,6,7,8}.
At the same time, the studies related to grey system forecasting models have been still emerged, and the prediction performance of grey forecasting models will be continuously enhanced as well \cite{9,10}.
In order to further improve the accuracy and application range of the grey models already in place, many scholars have proposed several new types of grey prediction models.
For example, Lang et al. \cite{11} proposed a new time-lag grey model to forecast photovoltaic power generation in the Asia-Pacific region. Zhou et al. \cite{12} put forward a new discrete grey model for natural gas prediction in Jiangsu Province, China. Wu et al. \cite{13}  constructed a new nonlinear grey model for short-term natural gas prediction. Zhou et al. \cite{14} established a novel discrete grey seasonal model and used several practical examples to verify the effectiveness of their model. Wang and Jv \cite{15} combined the grey prediction model and quantile regression so as to design a combination for forecasting time-series sequence. Xiao et al. \cite{16}  proposed an optimized nonlinear grey model for biomass energy consumption prediction.
 Li et al. \cite{17} searched the optimal fractional accumulation order for improved seasonal grey model by virtue of the particle swarm optimization (PSO). Zeng et al. \cite{18} optimized the structure of the grey Verhulst model.
Wei and Xie\cite{19}  proposed a method family for nonlinear parameter estimation of grey forecasting models, their numerical results show that the proposed method has advantages of high accuracy and robustness.
Zeng \cite{20}  proposed a novel discrete grey Riccati model and verified its effectiveness in a range of practical application scenarios.
Xiao et al. \cite{21}  put forward a grey Riccati Bernoulli model and successfully applied it to deal with clean energy consumption.
Ye et al. \cite{7} proposed a novel time-delay multivariate grey model
And used in China's carbon dioxide analysis.

Admittedly, the above models achieve satisfactory results and address many practical problems, however, the accumulation order in these models are all integer-order, which is prone to impair the prediction performance. In other words, the accumulation order should dynamically alter in accordance with modeing situations.
To this end, Wu et al. \cite{10} proposed the fractional gray model in 2013. Subsequently, various extensive forms of fractional models were proposed \cite{22,23}. At the same time, several other types of fractional accumulation operators have also been proposed, each with its own advantages. Ma et al. \cite{24}  proposed a conformable fractional accumulation (denoted as CFA) and built a conformable fractional grey prediction model. Chen et al. \cite{25} incorporated Fractional Hausdorff derivative into the modeling procedure of grey model so as to construct a new method. Liu et al. \cite{26} proposed the damping accumulated generating operator and constructed a new fractional grey model. With the continuous development of fractional grey models, fractional derivatives are applied on the differential equation of grey models \cite{27,28} .

For example, Mao et al. \cite{29} proposed a fractional grey model based on non-singular exponential kernel.\begin{equation}
\frac{1}{{1 - r}}\int_0^t {{x^\prime }} (\tau )\exp \left[ { - \frac{{r(t - \tau )}}{{1 - r}}} \right]d\tau  + a{x^{(r)}}(t) = b;{\text{   }}r > 0,t > 0.
\end{equation}
where ${x^{(r)}}(t)$ is the $r$-order accumulation \cite{10} of the original sequence ${x^{({\text{0}})}}(t)$,and literature \cite{28} proposed a fractional prediction model with conformable derivatives. The research status of fractional grey model can be seen in  Table \ref{Schoolar}. Using different fractional derivatives in the gray model can make full use of different types of derivatives to dig out potential laws among data. Fractional calculus has become a research hotspot of scholars \cite{29}. The famous Riemann-Liouville fractional integral with $r \in \mathbf{R}$ order \cite{30} by
\begin{equation}
I_{{t_0},r}^{RL}f(t) = \frac{1}{{\Gamma (r)}}\int_{t_0}^t {{{(t - \tau )}^{r  - 1}}} f(\tau )d\tau= {\Phi _r }*f(t);{\text{   }}r > 0,x > {t_0},
\end{equation}
motivated by the Cauchy integral
\begin{equation}
\int_{t_0}^t d {\tau _1}\int_{t_0}^{{\tau _1}} d {\tau _2} \cdots \int_{t_0}^{{\tau _{n - 1}}} f \left( {{\tau _n}} \right)d{\tau _n} = \frac{1}{{\Gamma (n)}}\int_{t_0}^t {{{(t - \tau )}^{n - 1}}} f(\tau )d\tau.
\end{equation}
Where ${\text{*}}$ is the convolution operator and ${\Phi _r}$ is the power function ${\Phi _r}(\tau) = \frac{1}{{\Gamma (r)}}\tau_ + ^{r - 1}$. Zhao and Luo \cite{30} proposed a generalized fractional operator, which proved that several types of calculus expressions are special cases of this operator.

\begin{table}[htbp]\footnotesize
  \centering
  \caption{Research findings of the different fractional grey models in recent years.}
    \begin{tabular}{lp{16.815em}lll}
    \toprule
    Author's name(year)   & \multicolumn{1}{l}{Approach} & Abbreviation &  FAGO & \multicolumn{1}{p{4.19em}}{FD} \\
    \midrule
    Wu (2013)\cite{10} & Grey model with the fractional\newline{}order accumulation & FGN(1,1) & \checkmark & - \\
    Mao (2016)\cite{28}  & Fractional grey system model & FGM(q,1)  & \checkmark & \checkmark \\
    Yang (2016)\cite{11} & Continuous fractional-order\newline{}grey model & GM(q,1)/GM(q,N) & -     & \checkmark \\
    Yang (2018)\cite{31} & \multicolumn{1}{l}{Interval grey modelling based on fractional calculus} & Interval GM(1,1)  &       & - \\
    Wu (2018)\cite{32} & The GMC(1,n) model with\newline{}fractional order accumulation & FGMC(1,n) & \checkmark & - \\
    Wu (2018)\cite{33} & Fractional order accumulation\newline{}grey model &FAGMO(1,1,k) & \checkmark & - \\
    Ma (2019)\cite{34} &  Fractional discrete multivariate\newline{}grey model & FDGM(1,n) & \checkmark & - \\
    Chen (2019)\cite{35} &  Time-delayed polynomial\newline{}fractional order grey model & TDPFOGM(1,1) & \checkmark & - \\
    Ma (2019)\cite{24} & \multicolumn{1}{l}{Conformable fractional grey model} & CFGM(1,1) & \checkmark & - \\
    Wang (2019)\cite{36} &  Fractional calculus function of\newline{}grey prediction model & FGM(1,1) & \checkmark & - \\
    Ma (2019)\cite{37} &  Fractional time delayed grey\newline{}model & FTDGM(1,1) & \checkmark & - \\
    Zhu (2019)\cite{38} &  Self-adaptive fractional\newline{}weighted grey model & SFOGM(1,1) & -     & - \\
    Yan (2020)\cite{25} & Fractional Hausdorff grey model & FHGM(1,1) & \checkmark & - \\
    Xie (2020)\cite{39} & Continuous grey model with conformable fractional derivative & CCFGM(1,1) & \checkmark & \checkmark \\
    Wu (2020)\cite{13} & conformable fractional non-homogeneous grey model & CFNGM(1,1) & \checkmark & - \\
    Kang(2020)\cite{40} & Fractional derivative multivariable grey model  & CFGMC(q,N) & \checkmark & \checkmark \\
    Gao(2020) \cite{41} &  fractional grey Riccati model & FGRM(1,1) & \checkmark & - \\
    Mao (2020)\cite{29}  & Fractional grey model based on non-singular exponential kernel & EFGM(q,1) & \checkmark & \checkmark \\
    Liu (2021)\cite{26} & The damping accumulated grey model & DAGM  & \checkmark & - \\
    Liu (2021)\cite{42} & Optimized fractional grey model-based variable background value & OFAGM(1,1) & \checkmark & - \\
    Zhang (2021)\cite{43} & Fractal derivative fractional grey Riccati model & FDFGRM & \checkmark & \checkmark \\
    \bottomrule
    \end{tabular}%
  \label{Schoolar}%
\end{table}%

Based on the above acknowledge, it is concluded that the fractional grey model has been widely used in diverse fields and a large number of fractional grey models has been constantly emerged. However, in return, there exist a variety of fractional grey models with various constructions of fractional accumulated generating operator, which is not conducive for beginners to understand and grasp the relevant modeling mechanism; at the same time, it is difficult for participants to identify which model can provide satisfactory results in specific scenarios, thus restraining the further methodological developments and practical applications.

To this end, this paper develops an unified framework for commonly-used fractional grey models by combining general fractional derivative with memory effects and grey modeling theory. The main contributions can be outlined as follows.
(i) In accordance with the GC and GRL derivatives, the unified framework for diverse fractional grey model is developed, this framework also derive other new fractional grey models by taking different kernel functions and normalization functions.
(ii) The exact solution (often called time response function) and the corresponding estimates of the structural parameters are, in detail, deduced from this proposed framework.
(iii) The four popular intelligent algorithms are employed to determine the emerging coefficients $r$, $\alpha$. (iv) Two published cases are utilized to verify the effectiveness of the UFGM(1,1) model and especially explore the effects of fractional accumulation order and initial value on the prediction performance. (v) The UFGM(1,1) model is also applied in two real cases to explain how to apply the proposed approach.

The rest of paper is organized as follows: Section \ref{sec:2} briefly introduces general fractional derivatives with memory effects. Section \ref{sec:3} discusses the modeling procedure of the unified framework. Section \ref{sec:4} verifies the validity of the UFGM(1,1) model. Section \ref{sec:5} applies this method in dealing with Henan's and Chongqing's water supply production capacity and Section \ref{sec:6} concludes.

\section{General fractional derivatives with memory effects}
\label{sec:2}
In this section, we will give a brief review of general fractional derivatives with memory effects and, in detail, deduce two important properties for these two kinds of general fractional derivative, which is helpful in constructing the unified framework for general fractional grey models in the following section.

\begin{definition}[See \cite{30}]
Given order $ r \in (0,1) , t>{t_0} $, the general fractional derivative with memory effects in the Caputo sense is presented as
\begin{equation}
D_{{t_0},r}^{GC}f(t) = N(r)\int_{t_0}^t k (t - \tau ,r)\frac{{df(\tau )}}{{d\tau }}d\tau,
\label{equ:GCD}
\end{equation}
where $G$ and $C$ represent the general definition and Caputo sense, respectively. $k(t - \tau ,r)$ and  $N(r)$ denote kernel function and normalization one, respectively.
\label{ref:1}
\end{definition}

Eq.(\ref{equ:GCD}) can be rewritten in the form of convolution as follows,
\begin{equation}
D_{{t_0},r}^{GC}f(t) = N(r)k(\tau ,r)*\frac{{df(t)}}{{dt}}
\label{equ:GCD_Cov}
\end{equation}

By taking different kernel functions and normalization ones, Eq.(\ref{equ:GCD_Cov}) can be transferred into diverse fractional derivatives, for example
\begin{eqnarray}\label{equ:Caputo}%
D_{{t_0},r}^{{{Caputo}}}f(t) &=& \frac{1}{{\Gamma (1 - r)}}\int_{t_0}^t {{{(t - \tau )}^{ - r}}} \frac{{df(\tau )}}{{d\tau }}d\tau\nonumber\\
&=&\frac{1}{{\Gamma (1 - r)}}{{\text{t}}^{ - r}}{\text{*}}\frac{{df(t)}}{{dt}}
\end{eqnarray}

Eq.(\ref{equ:Caputo}) is the basic form of the Caputo fractional derivative.

\begin{definition}[See \cite{30}]
Given $ r \in (0,1), t>{t_0} $, the general fractional derivative with memory effects in the Riemann$-$Liouville sense is represented as
\begin{equation}
D_{{t_0},r}^{GRL}f(t) = \frac{d}{{dt}}N(r)\int_{t_0}^t k (t - \tau ,r)f(\tau )d\tau,
\label{eq:grl}
\end{equation}
where $RL$ denotes the Riemann-Liouville sense, $k(t - \tau ,r)$ and  $N(r)$ are the same as in Definition \ref{ref:1}.
\label{ref:2}
\end{definition}

Similar to Eq.(\ref{equ:GCD_Cov}), Eq.(\ref{eq:grl}) can be expressed in the form of convolution as follows,
\begin{equation}
D_{{t_0},r}^{GRL}f(t) =\frac{d}{{dt}}N(r)k(\tau ,r)*f(t)
\label{eq:grl1}
\end{equation}

Meanwhile, Eq.(\ref{eq:grl1}) can be reduced to other fractional derivatives by taking different kenel functions and normalization ones, for example, the
Riemann-Liouville fractional derivative is calculated as
\begin{eqnarray}
D_{{t_0},r}^{{{RL }}}f(t) &=& \frac{{df(t)}}{{dt}}\frac{1}{{\Gamma (1 - r)}}\int_{t_0}^t {{{(t - \tau )}^{ - r}}} f(\tau )d\tau \nonumber \\
&=& \frac{{df(t)}}{{dt}}\frac{1}{{\Gamma (1 - r)}}{t^{ - r*}}f(t)
\label{equ:Caputo}%
\end{eqnarray}

Each fractional derivative has its own advantages and limitations, in addition, other kinds of the fractional derivative with diverse kernel functions and normalization ones are also shown in Ref.\cite{30} and references therein.

\section{Unified representation of fractional grey model}
\label{sec:3}
In this section, we present a unified framework for the general fractional grey models based on the above findings. In addition to analyze the characteristics of this framework, another focus is to deduce its time response function and structural parameters.

Assume $X^{(0)}=\left\{x^{(0)}(1),x^{(0)}(2),\cdots,x^{(0)}(n)\right\},n\geq4$ to be a nonnegative sequence, and the $\alpha$-order accumulated generating operation ($\alpha-FAGO$) sequence is given as $X^{(\alpha)}=\left\{x^{(\alpha)}(1),x^{(\alpha)}(2),\cdots,x^{(\alpha)}(n)\right\}$. then consider the general fractional grey model in the Caputo and Riemann-Liouville senses, namely UFGM, one has
\begin{equation}
N(r)\int_{t_0}^t k (t - \tau ,r)\frac{{d{x^{(\alpha )}}(\tau )}}{{d\tau }}d\tau  + a{x^{(\alpha )}}(t) = \psi(t) + c,
\label{eq:7}
\end{equation}
and
\begin{equation}
\frac{d}{{dt}}N(r)\int_{t_0}^t k (t - \tau ,r){x^{(\alpha )}}(\tau )d\tau  + a{x^{(\alpha )}}(t) = \psi(t) + c
\label{eq:8}
\end{equation}
respectively, where $\psi(t)$ is the unknown function with respect to $t$ and the discrete form of  ${x^{(\alpha )}}(t)$  gives
\begin{equation}\label{eq:fago}
{x^{(\alpha )}}(k) = \sum\limits_{i = 1}^k \binom{k - i + \alpha  - 1}{k-i} {x^{(0)}}(i),k = 1,2, \ldots ,n,
\end{equation}
where $\binom{k - i + \alpha  - 1}{k-i} =\frac{(\alpha+k-i-1)(\alpha+k-i-2) \cdots(\alpha+1) \alpha}{(k-i) !}$. The fractional accumulation plays an important role in the construction of grey models. For more details on the accuracy analysis of fractional accumulation, see \cite{10} for a review.

For simplicity, $\psi(t)$ is taken as $\psi(t) = bt + c$ and, Eqs.(\ref{eq:7})-(\ref{eq:8}) become
\begin{equation}
N(r) \int_{t_0}^t k (t - \tau ,r)\frac{{dx^{(\alpha)}(\tau)}}{{d\tau }}d\tau  + a{x^{(\alpha )}}(t) = bt + c,
\label{Raw_GC}
\end{equation}
and
\begin{equation} \label{GRL_grey}
\frac{d}{{dt}} N(r) \int_{t_0}^t k (t - \tau ,r)x^{(\alpha)}(\tau)d\tau  + a{x^{(\alpha )}}(t) = bt + c.
\end{equation}

Next, using the trapezoid formula gives the discrete-time equation of Eq.(\ref{Raw_GC})
\begin{equation}
\nabla D_{{t_0},r}^{GC}{x^{(\alpha )}}({\text{k}}) + \frac{a}{2}\left\{ {{x^{(\alpha )}}(k) + {x^{(\alpha )}}(k - 1)} \right\} = bk + c,k = 1,2,3,...n
\end{equation}

For the structural parameter estimation, set $\bm{\varphi} = {\left( {a,b,c} \right)^{\text{T}}}$ and $f\left(x^{(\alpha)}(k),\bm{\varphi} \right) =  - \frac{a}{2}\left( x^{(\alpha )}(k) + x^{(\alpha )}(k - 1) \right) + bk + c$, one has
\begin{eqnarray}
\bm{\varphi} = {\left( {a,b,c} \right)^{\text{T}}}=\left(\bm{A} ^{\text{T}}\bm{A}\right)^{-1}\bm{A} ^{\text{T}}\bm{Y}
\end{eqnarray}
where
\begin{equation}
A = \left[ {\begin{array}{*{20}{c}}
  {\frac{{\partial f\left( {{x^{(\alpha )}}(2),\varphi } \right)}}{{\partial a}}}&{\frac{{\partial f\left( {{x^{(\alpha )}}(2),\varphi } \right)}}{{\partial b}}}&{\frac{{\partial f\left( {{x^{(\alpha )}}(2),\varphi } \right)}}{{\partial c}}} \\
  {\frac{{\partial f\left( {{x^{(\alpha )}}(3),\varphi } \right)}}{{\partial a}})}&{\frac{{\partial f\left( {{x^{(\alpha )}}(3),\varphi } \right)}}{{\partial b}}}&{\frac{{\partial f\left( {{x^{(\alpha )}}(3),\varphi } \right)}}{{\partial c}}} \\
   \vdots & \vdots & \vdots  \\
  {\frac{{\partial f\left( {{x^{(\alpha )}}(n),\varphi } \right)}}{{\partial a}}}&{\frac{{\partial f\left( {{x^{(\alpha )}}(n),\varphi } \right)}}{{\partial b}}}&{\frac{{\partial f\left( {{x^{(\alpha )}}(n),\varphi } \right)}}{{\partial c}}}
\end{array}} \right],Y = \left[ {\begin{array}{*{20}{c}}
  {\nabla D_{{t_0},r}^{GC}{x^{(\alpha )}}(2)} \\
  {\nabla D_{{t_0},r}^{GC}{x^{(\alpha )}}(3)} \\
   \vdots  \\
  {\nabla D_{{t_0},r}^{GC}{x^{(\alpha )}}(n)}
\end{array}} \right],
\end{equation}
and $\nabla D_{{t_0},r}^{GC}{x^{(\alpha )}}(k)$ represents the discrete form of $D_{{t_0},r}^{GC}{x^{(\alpha )}}(k)$. It is worth noting that the structural parameter estimation of Eq.(\ref{GRL_grey}) has been omitted because of their high similarity.

\begin{theorem}
Set ${\text{ }}{\mathcal{L}^{ - 1}}\left\{ {\frac{1}{{\left[ {N(r)K(s,r)s+a} \right]}}} \right\} = \Theta (r,t,a)$ and $j$ indicates the discrete variable corresponding to $t$, the exact solution  of Eq.(\ref{Raw_GC}) is expressed as
\begin{equation} \label{GC_response}
{{\hat x}^{(\alpha)}}(j) = \Theta (r,j,a)*jb + \Theta (r,j,a)*c+{{\hat x}^{(\alpha)}}(0)N(r)k(j,r)*\Theta (r,j,a),
\end{equation}
where
\begin{equation}
{{\hat x}^{(\alpha)}}(0)=\frac{{{x^{(\alpha)}}(1) - \Theta (r,{\text{1}},a)*b + \Theta (r,{\text{1}},a)*c}}{{N(r)k({\text{1}},r)*\Theta (r,{\text{1}},a)}},
\end{equation}
 ${\mathcal{L}^{ - 1}}\left\{ . \right\}$ denotes the inverse Laplace transform. ${{\hat x}^{(r)}}(j)$ is the output value by the prediction model.
\end{theorem}

{\it\textbf{Proof.}} By applying the Laplace transform of general fractional derivative mentioned in \cite{30}, one has
\begin{equation} \label{Lapgrey}
N(r)\left\{sX(s) - {x^{(\alpha)}}(0)\right\}K(s,r ){\text{ + a}}X(s) = \frac{b}{{{s^2}}} + \frac{c}{s},
\end{equation}

Let $\frac{1}{{\left\{ {N(r)K(s,r)s{\text{ + }}a} \right\}}} = \Xi (r,s,a)$, Eq.(\ref{Lapgrey}) becomes
\begin{equation}\label{eq:18}
X(s) = \frac{b}{{{s^2}}}\Xi (r,s,a) + \frac{c}{s}\Xi (r,s,a) + {x^{(\alpha )}}(0)K(s,r)N(r)\Xi (r,s,a),
\end{equation}

Subsequently, through the inverse Laplace transform, Eq.(\ref{eq:18}) yields
\begin{eqnarray}\label{eq:19}
  {x^{(\alpha )}}(t) &=& {\mathcal{L}^{ - 1}}\{ \Xi (r,s,a)\} *tb + {\mathcal{L}^{ - 1}}{\{ \Xi (r,s,a)\} ^*}c \nonumber \\
  &+& {x^{(\alpha)}}(0)N(r){\mathcal{L}^{ - 1}}{\{ K(s,r)\} ^*}{\mathcal{L}^{ - 1}}\{ \Xi (r,s,a)\}
\end{eqnarray}
where $\mathcal{L}\left\{ {{{\hat x}^{(\alpha)}}(t)} \right\} = X(s)$. Then substituting Eq.(\ref{eq:19}) into Eq.(\ref{Lapgrey}), one can write
\begin{equation}
{{\hat x}^{(\alpha)}}(t) = \Theta (r,t,a)*tb + \Theta (r,t,a)*c{\text{ + }}{{\hat x}^{(\alpha)}}(0)N(r)k(t,r)*\Theta (r,t,a).
\end{equation}
Set ${{\hat x}^{(\alpha)}}(1) = {x^{(\alpha)}}(1)$ and $j=t$, this completes the proof. \hfill $\Box$

\begin{theorem}
The exact solution of Eq.(\ref{GRL_grey}) is expressed as
\begin{equation}
{\hat x^{(\alpha )}}(t) = \Theta(r,t,a)*tb + \Theta{(r,t,a)^*}c + {x^{(\alpha)}}(0)N(r)k{(x,r)^*}\Theta(r,t,a)
\label{GRLT}
\end{equation}
where $\widehat{x}^{(\alpha )}(0)=\frac{x^{(\alpha )}(1)-\Theta(r,1,a)*b + \Theta(r,1,a)^*c}{N(r)k(x,r)^*\Theta(r,1,a)}$ and $\mathcal{L}^{ - 1}\left\{ \frac{1}{N(r)K(s,r)s + x^{(\alpha )}(0) + a}\right\}  = \Theta(r, t, a)$.
\end{theorem}

{\it\textbf{Proof.}} By applying the link between the GC and GRL derivatives in \cite{30}
\begin{equation}
D_{r,{t_0}} ^{GC}f(t) = D_{r,{t_0}} ^{GRL}x^{(\alpha)}(t) - N(r )x^{(\alpha)}(0)k(x,r )
\end{equation}

Thus, Eq.(\ref{GRL_grey}) is rearranged as
\begin{equation}
D_{r,{t_0}}^{GC}x^{(\alpha)}(t)+N(r)x^{(\alpha)}(0)k(x,r)+ax^{(\alpha )}(t)=bt + c
\label{GRLmm}
\end{equation}

By applying the laplace transform on both sides of Eq. (\ref{GRLmm}), we obtain
\begin{equation} \label{GRLlaplace}
N(r)\left\{ sX(s) - x^{(\alpha )}(0) \right\}K(s,r) + x^{(\alpha )}(0)X(s) + aX(s) = \frac{b}{s^2} + \frac{c}{s}
\end{equation}

Set $\frac{1}{{\left\{ {N(r)K(s,r)s + {x^{(\alpha )}}(0) + a} \right\}}} = \mathcal{W}(r,s,a)$, it follows that
\begin{equation}
X(s) = \mathcal{W}(r,s,a)\frac{b}{{{s^2}}} + \mathcal{W}(r,s,a)\frac{c}{s} + {x^{(\alpha )}}(0)N(r)K(s,r)\mathcal{W}(r,s,a)
\end{equation}
and
\begin{eqnarray}
  \widehat {x}^{(\alpha )}(t) &=& \mathcal{L}^{ - 1}\left\{\mathcal{ W}(r,s,a)\right\} *tb + \mathcal{L}^{ - 1}\left\{\mathcal{ W}(r,s,a)\right\} *c \nonumber\\
&+& x^{(r)}(0)N(r)\mathcal{L}^{-1}\left\{K(s,r)\right\}^*\mathcal{L}^{-1}\left\{\mathcal{W}(r,s,a)\right\}
\end{eqnarray}
where $ \mathcal{L}^{ - 1}\left\{ \mathcal{W}(r,s,a)\right\}  = \Theta(r,t,a)$. Finally, set  ${\hat x^{(\alpha )}}(1) = {x^{(\alpha )}}(1)$ and $s=t$, the exact solution (often called time response function) of Eq.(\ref{GRLT}) can be easily obtained. \hfill $\Box$

\begin{remark}\label{re1}
Set $k(t-\tau,r)=e^{-\frac{r}{1-r}x}$ and $N(r)=1$ \cite{29}, a novel fractional grey model is presented as
\begin{equation} \label{grey_exp}
_0^{\exp }D_t^r{x^{(\alpha )}}(t) + a{x^{(\alpha )}}(t) = bt + c.
\end{equation}
\end{remark}

By applying the laplace transformation on Eq.(\ref{grey_exp}), which yields that
$$
\mathcal{L}\left\{ {_0^{\exp }D_t^r{x^{(\alpha )}}} \right\} + \mathcal{L}\left\{ {a{x^{(\alpha )}}(t)} \right\} = \mathcal{L}\left\{ {bt} \right\} + \mathcal{L}\left\{ c \right\}
$$
and
\begin{eqnarray}\label{eq:las}
\frac{{sX(s) - {x^{(\alpha )}}(0)}}{{s + r(1 - s)}} + {a}X(s) = \frac{{{b}}}{{{s^2}}} + \frac{{{c}}}{s}
\end{eqnarray}

By further solving Eq.(\ref{eq:las}), one can write
\begin{equation}\label{eq:xs}
X(s) = \frac{{cr{s^2} + brs - crs - c{s^2} - {x^{(\alpha )}}(0){s^2} - br - bs}}{{{s^2}\left( {ars - ar - as - s} \right)}}
\end{equation}

The continuous solution to Eq.(\ref{eq:xs}) can be calculated by the inverse laplace transformation operation, which is as follows,
\begin{equation}
{\widehat x^{(\alpha )}}({\text{t}}) = \frac{{b\left( {\frac{{art}}{{ar - a - 1}} - 1} \right)}}{{{a^2}r}} + \frac{{bt(ar - a - 1) + {e^{\frac{{art}}{{ar - a - 1}}}}( - ax + c) + c(ar - a - 1)}}{{(ar - a - 1)a}}
\end{equation}

Set ${\widehat x^{(\alpha )}}(1) = {x^{(\alpha )}}(1)$ and $t=k$, the discrete time response function of \label{grey_exp} can be given as
\begin{equation} \label{exp_trsponse_contiue}
{\hat x^{(\alpha )}}(k) = \left( {{x^{(\alpha )}}(1) - \frac{b}{a} - \frac{c}{a} + \frac{b}{{{a^2}r}}} \right){e^{\frac{{ - ar(k - 1)}}{{1 - ar + a}}}} + \frac{b}{a}k + \frac{c}{a} - \frac{b}{{{a^2}r}}
\end{equation}

\begin{remark}\label{re2}
Inspired by the idea in \cite{28},
the discrete forma of Eq.(\ref{grey_exp}) is
$
{x^{(\alpha  - r)}}(k) + \frac{a}{2}\left\{ {{x^{(\alpha )}}(k - 1) + {x^{(\alpha )}}(k)} \right\}= bk + c
$.
Subsequently we estimate the structural parameters by constructing a cost function expressed as
\begin{equation}\label{eq:of}
{\chi ^2} = {\sum\limits_{k = 2}^n {\left\{ {{x^{(\alpha  - r)}}(k) + a{z^{(\alpha )}}(k) - bk - c} \right\}} ^2}
\end{equation}
\end{remark}

And by minimizing ${\chi ^2}$, one has a set of equations with respect to partial function expressed as
\begin{equation} \label{para_estimation}
\left\{ {\begin{array}{*{20}{l}}
  {\frac{{\partial {\chi ^2}}}{{\partial a}} = \sum\limits_{k = 2}^n {2{z^{(\alpha )}}(k){x^{(\alpha  - r)}}(k) + 2a{{\left\{ {{z^{(\alpha )}}(k)} \right\}}^2} - 2{z^{(\alpha )}}(k)bk - 2{z^{(\alpha )}}(k)c}  = 0,} \\
  \begin{gathered}
  \frac{{\partial {\chi ^2}}}{{\partial b}} = \sum\limits_{k = 2}^n { - 2k{x^{(\alpha  - r)}}(k) - 2ka{z^{(\alpha )}}(k) + 2{k^2}b + 2kc}  = 0, \hfill \\
  \frac{{\partial {\chi ^2}}}{{\partial c}} = \sum\limits_{k = 2}^n { - 2{x^{(\alpha  - r)}}(k) - 2a{z^{(\alpha )}}(k) + 2bk + 2c = 0.}  \hfill \\
\end{gathered}
\end{array}} \right.
\end{equation}
which yields that
\begin{equation} \label{para2}
\left\{ {\begin{array}{*{20}{l}}
  {\mathop \sum \limits_{k = 2}^n a{{\left\{ {{z^{(\alpha )}}(k)} \right\}}^2} - \mathop \sum \limits_{k = 2}^n {z^{(\alpha )}}(k)bk - \mathop \sum \limits_{k = 2}^n {z^{(\alpha )}}(k)c = \mathop { - \sum }\limits_{k = 2}^n {z^{(\alpha )}}(k){x^{(\alpha  - r)}}(k),} \\
  \begin{gathered}
    - \mathop \sum \limits_{k = 2}^n ka{z^{(\alpha )}}(k) + \mathop \sum \limits_{k = 2}^n {k^2}b + \mathop \sum \limits_{k = 2}^n kc = \mathop \sum \limits_{k = 2}^n k{x^{(\alpha  - r)}}(k), \hfill \\
   - \mathop \sum \limits_{k = 2}^n a{z^{(\alpha )}}(k) + \mathop \sum \limits_{k = 2}^n bk + \mathop \sum \limits_{k = 2}^n c = \mathop \sum \limits_{k = 2}^n {x^{(\alpha  - r)}}(k). \hfill \\
\end{gathered}
\end{array}} \right.
\end{equation}
and then, using the least squares, the structural parameter estimates give in matrix form expressed as
\begin{eqnarray}
\left(a,b,c\right)^{\text{T}}=\left(\mathcal{O}^{\text{T}}\mathcal{O}\right)^{-1}\mathcal{O}^{\text{T}}{\rm H}
\end{eqnarray}
where
$$\mathcal{O}= \begin{pmatrix}
\sum_{k=2}^n\left\{z^{(\alpha)}(k)\right\}^2&-\sum_{k=2}^nz^{(\alpha)}(k)k&-\sum_{k=2}^nz^{(\alpha)}(k)\\
-\sum_{k=2}^nkz^{(\alpha)}(k)&\sum_{k=2}^nk^2&\sum_{k=2}^nk\\
-\sum_{k=2}^nz^{(\alpha)}(k)&\sum_{k=2}^nk&\sum_{k=2}^n
\end{pmatrix}$$

$$
{\rm H}=\begin{pmatrix}
 -\sum_{k=2}^nz^{(\alpha)}(k)x^{(\alpha-r)}(k)\\
  -\sum_{k=2}^nkx^{(\alpha-r)}(k)\\
   -\sum_{k=2}^nx^{(\alpha-r)}(k)
\end{pmatrix}
$$

Substituting the estimates of the structural parameters and initial value $\widehat x^{(\alpha)}(1)=x^{(0)}(1)$ into Eq.(\ref{exp_trsponse_contiue}), the restored values of the original sequence can be given by then using the inverse fractional accumulated generating operation
\begin{eqnarray}
 \label{RFA}
{\hat x^{(0)}}(k) = \sum\limits_{i = 1}^k {\left( {\begin{array}{*{20}{c}}
  {k - i + 1 - \alpha  - 1} \\
  {k - i}
\end{array}} \right)} {\hat x^{(\alpha )}}(i) - \sum\limits_{i = 1}^{k-1} {\left( {\begin{array}{*{20}{c}}
  {k - 1 - i + 1 - \alpha  - 1} \\
  {k - i - 1}
\end{array}} \right)} {\hat x^{(\alpha )}}(i)
\end{eqnarray}

We recall the above remarks, where Remark \ref{re1} shows that the unified model has the flexibility to present the continuous fractional grey models and integer-order grey models. For instance, if fractional accumulation order $\alpha$ is taken as ``1'', Eq.(\ref{exp_trsponse_contiue}) is simplified as the basic form of NHGM(1,1) \cite{3}, whose time response function is presented as ${\hat x^{(\alpha )}}(k) = \left( {{x^{(\alpha )}}(0) - \frac{b}{a} - \frac{c}{a} + \frac{b}{{{a^2}}}} \right){e^{ - a(k - 1)}} + \frac{b}{a}k + \frac{c}{a} - \frac{b}{{{a^2}}}$. Remark \ref{re2}, in detail, introduces the calculation process of the least squares in the special example, which help readers better understand the relevant calculation mechanism. Remarks \ref{re1}-\ref{re2} also introduce the possibility of developing some other novel fractional models.

\section{Simulation studies}
\label{sec:4}
For the purpose of numerically analyzing the prediction performance of the unified model (i.e., UFGM(1,1)), we carry out the experiments from the diverse two aspects: One is to compare the prediction results by UFGM(1,1) with those by other benchmarks including the GM(1,1) \cite{1}, DGM(1,1) \cite{3}, FGM(1,1) \cite{10} and ANN \cite{45} models; another is to explore the impact of initial value and fractional accumulation order on the prediction accuracy of UFGM(1,1).

\subsection{Validity of the UFGM(1,1) model}
\label{sec:4.1}
Three error-value metrics are chosen to evaluate the accuracy of UFGM(1,1), which embrace the mean absolute percentage error (MAPE), mean square error (MSE) and mean absolute error (MAE), which are defined, respectively, as follows,
\begin{eqnarray}
MAPE&=&\frac{1}{n}\sum_{k=1}^{n}\frac{\left|\widehat x^{(0)}(k)-x^{(0)}(k)\right|}{x^{(0)}(k)}\times 100\% \label{eq:mape}\\
MSE&=&\frac{1}{n}\sum_{k=1}^{n}\left(\widehat x^{(0)}(k)-x^{(0)}(k)\right)^2\\
MAE&=&\frac{1}{n}\sum_{k=1}^{n}\left|\widehat x^{(0)}(k)-x^{(0)}(k)\right|
\end{eqnarray}

\begin{figure}
\centering
\includegraphics[scale=0.3]{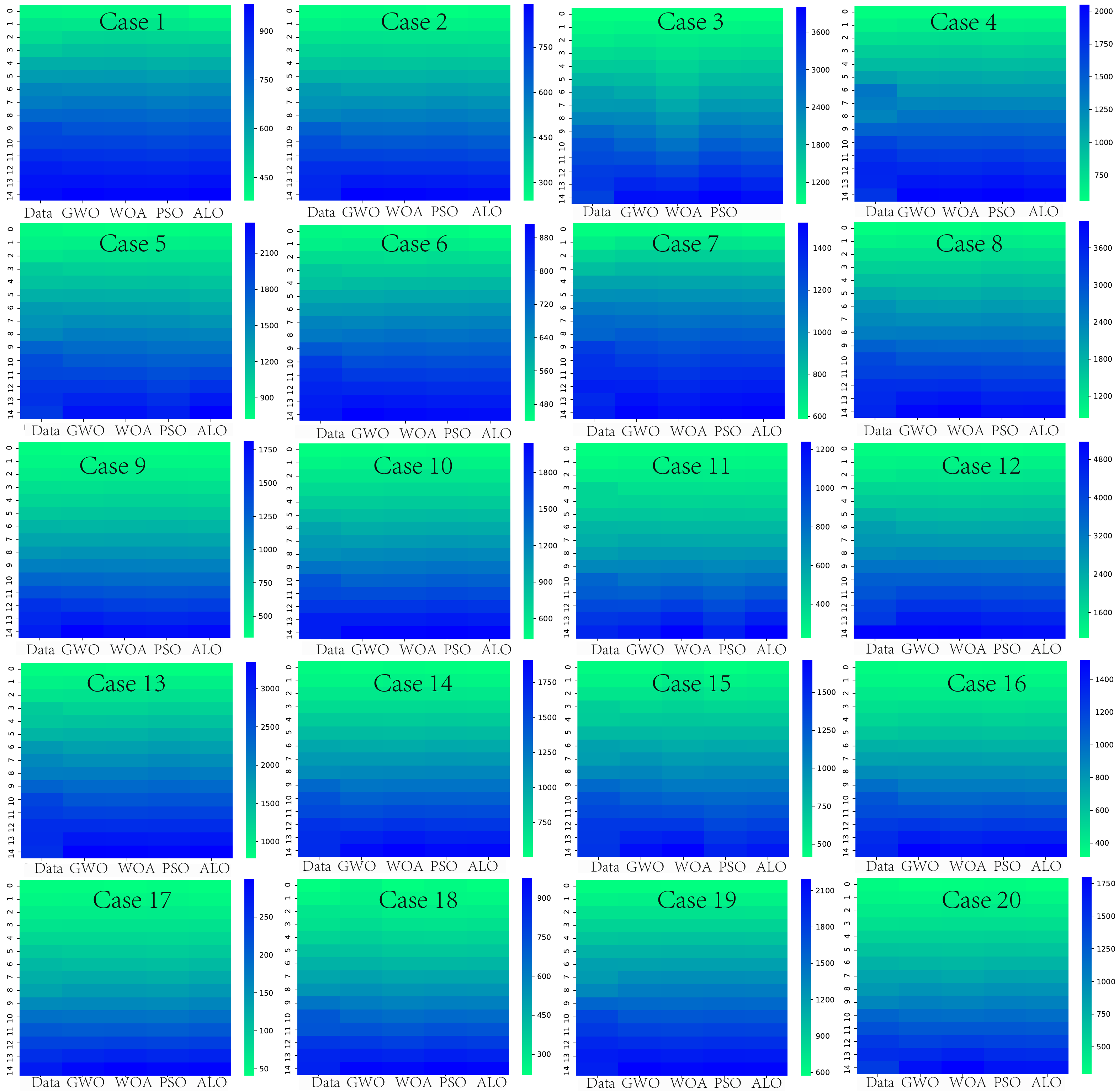}
\caption{Prediction performance of the different algorithms-based UFGM(1,1) model}
\label{fig:sear}
\end{figure}

For the experimental design, we consider 20 data sets, collected from the National Bureau of Statistics of China (http://www.stats.gov.cn/), as the inputs to calibrate models, which describes the electricity consumption in 20 provinces of China from 2014 to 2019. In UFGM(1,1), the emerging coefficients (i.e., $r$, $\alpha$) plays an important role in the modeling process and more importantly, has a directly impact on the prediction accuracy. Thus we utilize the four commonly-used intelligent techniques, namely the  grey wolf optimizer (GWO) \cite{46}, whale optimization algorithm (WOA)  \cite{47}, particle swarm optimization (PSO) \cite{48} and ant lion optimization (ALO) \cite{49}, to determine the optimal emerging coefficients for UFGM(1,1). For intuition, the fitting results produced by UFGM(1,1) with the help of the different algorithms are plotted in Figure \ref{fig:sear}.

Figure \ref{fig:sear} shows that there is no fixed algorithm that enables UFGM(1,1) to produce the highest accuracy in all data sets, the reason for this situation is that each algorithm has its own advantages and limitations, thus it is possible to obtain better solutions for different problems. In the following analysis, it is necessary to evaluate the fitness value of these algorithms in each experiment, and the optimal emerging coefficients are chosen for model calibration. By reference with the current study and other published papers related to benchmarks, their comparative results are shown in Table \ref{table:perf}. As table shows, the MAPE and MAE values of UFGM(1,1) are all smaller than those of other competitors, although its MSE values are slightly higher than traditional fractional grey model in several experiments, inferring that the UFGM(1,1) model outperforms other benchmarks in terms of model's prediction performance.

\begin{table}\scriptsize
  \centering
  \caption{Errors by the five competitive models in 20 data sets.}
    \begin{tabular}{cccccccc}
    \hline
    Case&Region & Metrics & GM(1,1) & DGM(1,1) & ANN & FGM(1,1) & UFGM(1,1) \\
    \hline
    \multirow{3}*{1}&\multirow{3}*{Beijing} &  MAPE & 4.4274 & 4.4237 & 3.362543 & 1.6294 & \textbf{1.4879} \\
         & & MSE   & 1046  & 1047.3 & 826.1047 & \textbf{210.14} & 235.57 \\
         & & MAE   & 27.941 & 27.86 & 24.12434 & 11.186 & \textbf{10.614} \\
        \multirow{3}*{2}&\multirow{3}*{Tianjin} &  MAPE & 4.4274 & 6.1891 & 3.132174 & 3.8301 & \textbf{2.404} \\
         & & MSE   & 1046  & 1312.1 & 495.232 & \textbf{533.62} & 646.51 \\
         & & MAE   & 27.941 & 28.989 & 18.21007 & 20.381 & \textbf{15.52} \\
    \multirow{3}*{3}&\multirow{3}*{Hebei } &  MAPE & 10.76 & 10.854 & 4.81092 & 7.9273 & \textbf{4.0727} \\
         & & MSE   & 54708 & 54986 & 14763.99 & \textbf{29176} & 33992 \\
         & & MAE   & 188.23 & 188.5 & 108.4177 & 149.38 & \textbf{104.99} \\
   \multirow{3}*{4}&\multirow{3}*{Shanxi } &  MAPE & 10.596 & 10.669 & 6.432003 & 6.8308 & \textbf{4.5779} \\
         & & MSE   & 17969 & 18005 & 9898.418 & \textbf{8890.9} & 9452.1 \\
         & & MAE   & 113.37 & 113.58 & 87.11586 & 81.838 & \textbf{65.915} \\
   \multirow{3}*{5}&\multirow{3}*{Liaoning} &  MAPE & 5.7687 & 4.88467 & 3.135815 & 4.0297 & \textbf{3.0828} \\
         & & MSE   & 10219 & 8775.161 & 6730.438 & \textbf{4893.4} & 6243.1 \\
         & & MAE   & 77.149 & 75.33215 & 66.32129 & 59.219 & \textbf{48.428} \\
    \multirow{3}*{6}&\multirow{3}*{Heilongjiang} &  MAPE & 2.4237 & 2.4298 & 2.524523 & 2.054 & \textbf{1.6633} \\
        &  & MSE   & 510.6 & 510.81 & 511.9586 & 291.74 & \textbf{260.93} \\
         & & MAE   & 16.947 & 16.962 & 18.7526 & 14.031 & \textbf{12.103} \\
   \multirow{3}*{7}&\multirow{3}*{Shanghai} &  MAPE & 7.5002 & 7.5101 & 2.947381 & 2.9327 & \textbf{2.1168} \\
         & & MSE   & 6967.2 & 6971.7 & 1679.166 & 1439.8 & \textbf{999.77} \\
         & & MAE   & 72.783 & 72.7  & 33.71344 & 31.477 & \textbf{23.533} \\
    \multirow{3}*{8}&\multirow{3}*{Zhejiang} &  MAPE & 10.609 & 10.652 & 3.280667 & 5.1743 & \textbf{3.1869} \\
         & & MSE   & 47904 & 48199 & 9964.663 & \textbf{13238} & 14965 \\
         & & MAE   & 183.32 & 182.89 & 85.00789 & 100.01 & \textbf{84.197} \\
    \multirow{3}*{9}&\multirow{3}*{Anhui} &  MAPE & 5.6001 & 5.6303 & 5.0247 & 4.2216 & \textbf{2.125} \\
          && MSE   & 3372  & 3449  & 3805.964 & 2134.8 & \textbf{1828.6} \\
          && MAE   & 44.015 & 43.871 & 53.12803 & 35.507 & \textbf{24.597} \\
   \multirow{3}*{10}&\multirow{3}*{Fujian} &  MAPE & 6.5022 & 6.5479 & 4.896582 & 3.9373 & \textbf{2.0374} \\
         & & MSE   & 5766.2 & 5839.6 & 5161.624 & 2369.6 & \textbf{2138.1} \\
          && MAE   & 56.337 & 56.409 & 57.57798 & 41.634 & \textbf{27.153} \\
    \multirow{3}*{11}&\multirow{3}*{Jiangxi} &  MAPE & 6.1053 & 6.1617 & 4.425646 & 3.7039 & \textbf{2.7821} \\
          && MSE   & 1254.7 & 1281.8 & 1303.201 & \textbf{574.3} & 641.04 \\
          && MAE   & 27.619 & 27.768 & 28.03915 & 19.493 & \textbf{16.943} \\
    \multirow{3}*{12}&\multirow{3}*{Shandong} &  MAPE & 7.1293 & 7.1613 & 5.36149 & 3.0638 & \textbf{2.9972} \\
          && MSE   & 31808 & 32177 & 36099.57 & 16175 & \textbf{15496} \\
          && MAE   & 148.94 & 148.89 & 160.1942 & 91.306 & \textbf{90.343} \\
   \multirow{3}*{13}&\multirow{3}*{Henan} &  MAPE & 9.5074 & 9.5784 & 4.283333 & 5.791 & \textbf{3.7368} \\
          && MSE   & 34638 & 34792 & 9638.134 & \textbf{15935} & 21391 \\
         & & MAE   & 155.34 & 155.53 & 85.30263 & 104.06 & \textbf{86.514} \\
    \multirow{3}*{14}&\multirow{3}*{Hubei} &  MAPE & 5.9726 & 6.018 & 4.077524 & 4.4207 & \textbf{2.6483} \\
          && MSE   & 5783  & 5809.9 & 3683.539 & 3511.3 & \textbf{2765.1} \\
          && MAE   & 60.408 & 60.544 & 49.47484 & 47.658 & \textbf{34.901} \\
    \multirow{3}*{15}&\multirow{3}*{Hunan} &  MAPE & 7.1714 & 7.2392 & 5.419488 & 4.7514 & \textbf{3.7075} \\
         & & MSE   & 5898.6 & 5922.9 & 4106.629 & \textbf{3018.3} & 3426.1 \\
          && MAE   & 62.315 & 62.555 & 55.92601 & 45.21 & \textbf{41.214} \\
    \multirow{3}*{16}&\multirow{3}*{Guangxi} &  MAPE & 7.3815 & 7.4647 & 4.084182 & 5.8011 & \textbf{3.3829} \\
          && MSE   & 4063.5 & 4111  & 1579.356 & 2460.7 & \textbf{2418.8} \\
          && MAE   & 50.793 & 50.964 & 35.79186 & 41.435 & \textbf{30.92} \\
    \multirow{3}*{17}&\multirow{3}*{Hainan} &  MAPE & 6.6034 & 6.6737 & 4.084182 & 3.5861 & \textbf{2.3553} \\
          && MSE   & 65.205 & 69.117 & 1579.356 & \textbf{25.571} & 26.381 \\
          && MAE   & 6.4182 & 6.4535 & 35.79186 & 4.2047 & \textbf{3.6773} \\
    \multirow{3}*{18}&\multirow{3}*{Chongqing} &  MAPE & 4.9022 & 4.9783 & 3.78875 & 4.5603 & \textbf{2.9425} \\
          && MSE   & 1036.4 & 1048.1 & 726.5574 & 796.52 & \textbf{775.06} \\
          && MAE   & 22.673 & 22.88 & 22.4828 & 21.756 & \textbf{17.092} \\
    \multirow{3}*{19}&\multirow{3}*{Sichuan} &  MAPE & 5.9849 & 6.0382 & 4.378797 & 4.1841 & \textbf{2.6998} \\
          && MSE   & 8572.9 & 8615  & 6206.165 & 4955.6 & \textbf{3959.6} \\
          && MAE   & 71.844 & 72.127 & 65.22631 & 58.053 & \textbf{42.153} \\
    \multirow{3}*{20}&\multirow{3}*{Yunnan } &  MAPE & 8.4035 & 8.5588 & 5.935477 & 7.3737 & \textbf{3.5109} \\
          && MSE   & 7395.2 & 7500.6 & 6465.405 & 5041.9 & \textbf{4563.6} \\
          && MAE   & 59.009 & 59.653 & 58.08669 & 58.986 & \textbf{36.112} \\
          \hline
    \end{tabular}%
  \label{table:perf}%
\end{table}%

\subsection{Effects of fractional accumulation order and initial value on modelling}
\label{sec:4.2}
In order to explore the effects of fractional accumulation order and initial value on modelling, the initial value is sampled at every interval of 200000 in the range of [1, 800000], and then $r$, $\alpha$, $a$, $b$, $c$ are all controlled within [0, 1], of which the step is 0.2. Consider the data on total population at the end of the year ($10^4$ persons) from the National Bureau of Statistics of China (http://www.stats.gov.cn/). The pseudo-code for the relevant calculations are shown in Algorithm \ref{alg:1}.

\begin{algorithm}\label{alg:1}
  \caption{Effects of the fractional accumulation order and initial value on the prediction performance.}
  \KwIn{Raw data $X ^{( 0 )}=\left\{x^{(0)}(1), x^{(0)}(2), \ldots, x^{(0)}(n)\right\}$}
  \KwOut{Simulation value ${\hat x^{(0)}}(k)$}
  Initialize ${\text{ MAPE}}{{\text{ }}_{\min }} =  + {\text{ inf }}$
    \For{$Initial$ $value=1$ to $800000$ $Step=200000$}
     {
  \For{$\alpha=0.1$ to $1$ $Step=0.2$}
  {
    \For{$r=0.1$ to $1$ $Step=0.2$}
      {

     \For{$a=0.1$ to $1$ $Step=0.2$}
     {
      \For{$b=0.1$ to $1$ $Step=0.2$}
      {
       \For{$c=0.1$ to $1$ $Step=0.2$}
    {
      Calculating ${X^{(\alpha )}} = \left( {{x^{(\alpha )}}(1),{x^{(\alpha )}}(2), \ldots ,{x^{(\alpha )}}(n)} \right)$ using the Eq.(\ref{eq:fago})\;
       Calculating $\hat{x}^{(\alpha)}(k)$ using the Eq.(\ref{exp_trsponse_contiue})\;
      Calculate the fitting results $\hat{x}^{(0)}(k)$ using Eq.(\ref{RFA})\;
      Calculate the error of this model using the Eq.(\ref{eq:mape})\;
    }
    }
    }
    }
  }
    return $\hat{x}^{(0)}(k)$\;
  }
\end{algorithm}

\begin{figure}
\centering
\includegraphics[scale=0.4]{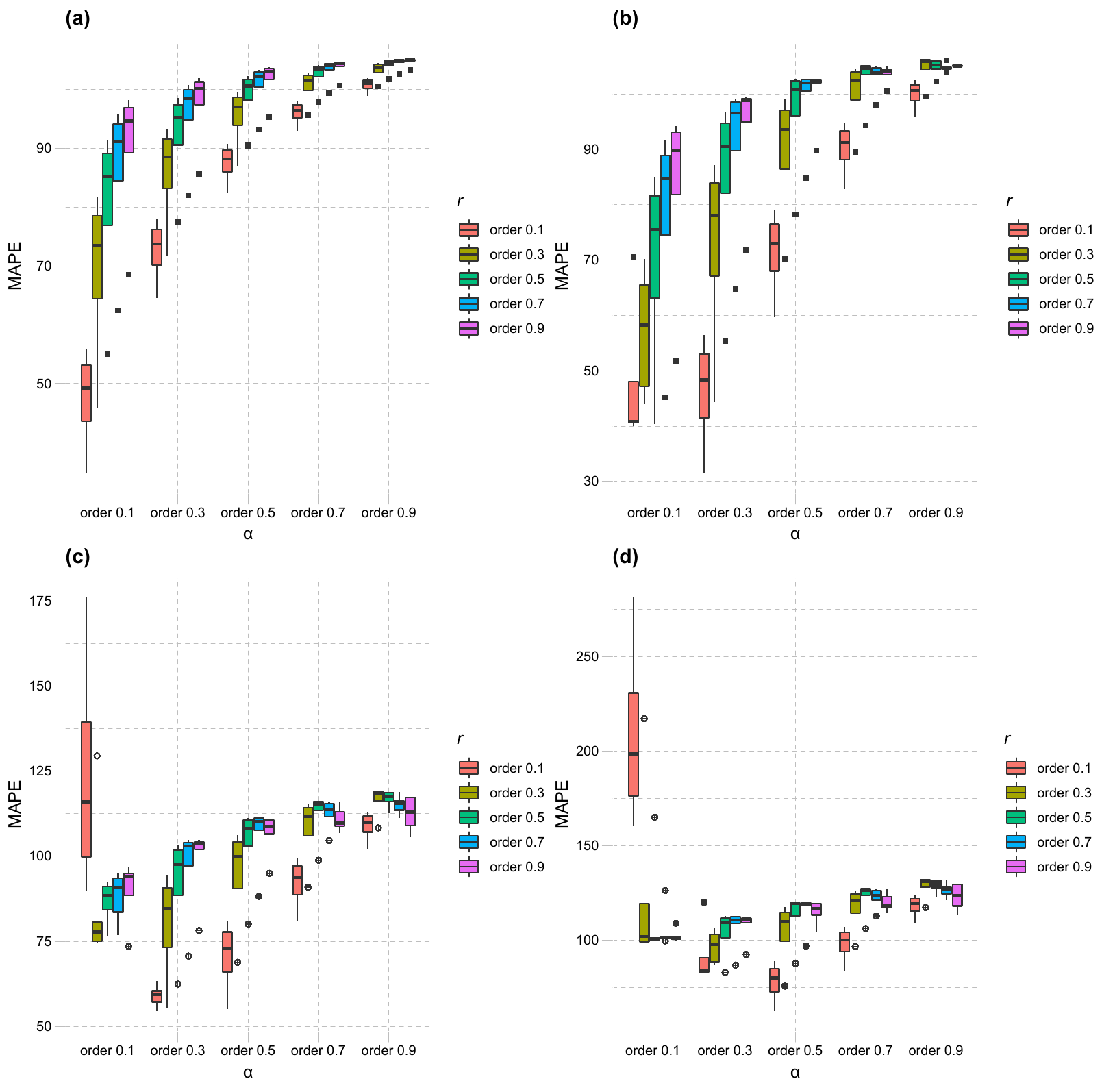}
\caption{Effects of the fractional order accumulation order $r,\alpha$ on the UFGM(1,1) model in the four initial value scenarios: (a)$initial\ value=1$;
 (b)$initial\ value=200001$; (c)$initial\ value=400001$; (d)$initial\ value=600001$;}
\label{fig:box}
\end{figure}

Figure \ref{fig:box} shows that with the increase of the initial condition value, the MAPE values get larger when the other conditions are fixed, for example, $r=0.1$ and $\alpha=0.1$, inferring the predicted values by the UFGM(1,1) model deviate from the actual ones. In return, in the small initial value scenario, the smaller the fractional accumulation order, the lower the fitting error (i.e., MAPE); with the increase of the initial value, the fractional accumulation order should properly increase so as to better mine the characteristics hidden in the original sequence, for example, $r=\alpha=0.3$ enables this model to provide the more accurate forecasts rather than that of $r=\alpha=0.1$ in the large initial value scenario (i.e., $initial\ value =800000$).

\section{Application}
\label{sec:5}
Foresight of water supply capacity is of great significance for the sustainable development of
society, is conducive to the rational distribution of water resources, and plays an important role in the
development of the national economy.
In this section, we further apply the UFGM(1,1) model to deal with the data on water supply production capacity (ten thousand cubic meters/day) in Henan and Chongqing, respectively. The prediction performance of the UFGM(1,1) model will be compared with the above benchmarks in terms of MAPE (elaborated on in Section \ref{sec:4.1}). The total amount of water supply production capacity is collected from  the National Bureau of Statistics of China (http://www.stats.gov.cn/), as shown in Tables \ref{table:c1}-\ref{table:c2}.

\begin{figure}
\flushleft
\includegraphics[scale=0.52]{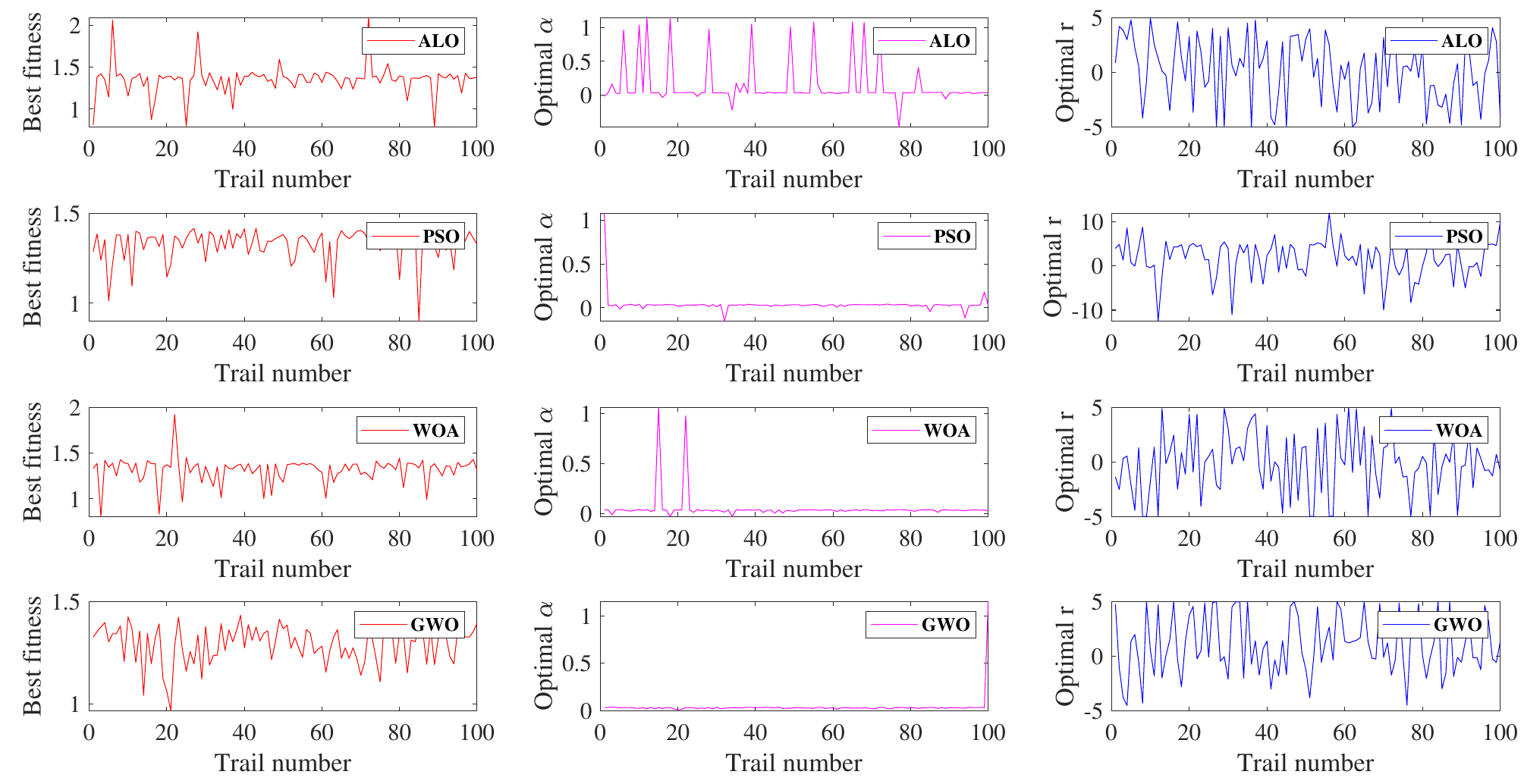}
\caption{$r$, $\alpha$ and the fitness value of the UFGM(1,1) model in each trial by using the four algorithms in Case 1.}
\label{fig:search}
\end{figure}

\textbf{Case 1.} (Forecasting Henan's water supply production capacity) We consider the total amount of water supply production capacity from 2004 to 2019, in which the data from 2004 to 2015 are used for model calibration, and the left four samples are used to examine the accuracy. For the purpose of obtaining the optimal emerging coefficients for the UFGM(1,1) model, the four algorithms are carried out in this experiments, thus the track of searching process by the different algorithms is graphed in Figure \ref{fig:search}. Subsequently, the minimum MAPE corresponding to the algorithm is applied to determine the emerging coefficients.

\begin{table}[htbp]\scriptsize
  \centering
  \caption{Simulated and predicted values by the five competitive models in Case 1.}
    \begin{tabular}{ccccccccccc}
    \hline
    Raw data & GM(1,1) & Error(\%) & DGM(1,1) & Error(\%) & ANN   & Error(\%) & FGM(1,1) & Error(\%) & UFGM(1,1) & Error(\%) \\
\hline
    \multicolumn{11}{c}{In-sample}\\
    \hline
    1038.31 & 1038.31 & 0.000     & 1003.8 & 3.32367 & 1038.31 & 0.0000     & 1038.31 & 0.000     & 1038.31 & 0.000 \\
    1026.51 & 1003.7 & 2.222092 & 1003.8 & 2.212351 & 1026.51 & 0.0000     & 1026.8 & 0.028251 & 1022.9 & 0.351677 \\
    1023.7 & 1011  & 1.240598 & 1011.1 & 1.230829 & 1023.7 & 0.0000     & 1021.8 & 0.185601 & 1014.9 & 0.859627 \\
    1039.85 & 1018.4 & 2.062798 & 1018.5 & 2.053181 & 1039.85 & 0.0000     & 1019.3 & 1.976247 & 1011.5 & 2.726355 \\
    1013.91 & 1025.9 & 1.182551 & 1026  & 1.192414 & 1047.962 & 3.358434 & 1018.7 & 0.472429 & 1011.9 & 0.198242 \\
    1007.79 & 1033.4 & 2.541204 & 1033.5 & 2.551127 & 1042.344 & 3.42871 & 1019.8 & 1.191717 & 1015.9 & 0.804731 \\
    1010.34 & 1041  & 3.034622 & 1041  & 3.034622 & 1039.884 & 2.924204 & 1023.1 & 1.262941 & 1023.3 & 1.282737 \\
    1037.56 & 1048.6 & 1.064035 & 1048.6 & 1.064035 & 1031.243 & 0.608813 & 1029.3 & 0.796099 & 1033.9 & 0.352751 \\
    1042.31 & 1056.3 & 1.342211 & 1056.3 & 1.342211 & 1034.192 & 0.778837 & 1039.8 & 0.240811 & 1047.6 & 0.507527 \\
    1047.26 & 1064.1 & 1.608006 & 1064  & 1.598457 & 1042.544 & 0.450299 & 1056.4 & 0.872754 & 1064.5 & 1.646201 \\
    1083.62 & 1071.9 & 1.08156 & 1071.8 & 1.090788 & 1055.567 & 2.588786 & 1082.1 & 0.140271 & 1084.3 & 0.062753 \\
    1121.39 & 1079.7 & 3.717707 & 1079.6 & 3.726625 & 1071.847 & 4.417999 & 1121.4 & 0.000892 & 1107.1 & 1.274311 \\
     MAPE     &       & 1.9183 &       & 1.9171 &       & 2.319511 &       & \textbf{0.65119} &       & 0.91541 \\
         \hline
    \multicolumn{11}{c}{Out-of-sample}\\
    \hline
    1180.32 & 1087.6 & 7.855497 & 1087.5 & 7.863969 & 1093.27 & 7.375093 & 1181.1 & 0.066084 & 1132.8 & 4.026027 \\
    1150.37 & 1095.6 & 4.761077 & 1095.5 & 4.76977 & 1132.02 & 1.595139 & 1271.4 & 10.52096 & 1161.3 & 0.950129 \\
    1166.64 & 1103.7 & 5.39498 & 1103.5 & 5.412124 & 1157.322 & 0.798721 & 1407.6 & 20.65419 & 1192.6 & 2.225194 \\
    1281.52 & 1111.7 & 13.25145 & 1111.5 & 13.26706 & 1181.88 & 7.77515 & 1612.6 & 25.83495 & 1226.7 & 4.277733 \\
    MAPE     &       & 7.8142 &       & 7.8292 &       & 4.386027 &       & 14.26904 &       & \textbf{2.8704} \\
\hline
    \end{tabular}%
  \label{table:c1}%
\end{table}%

\begin{figure}
\centering
\includegraphics[scale=0.6]{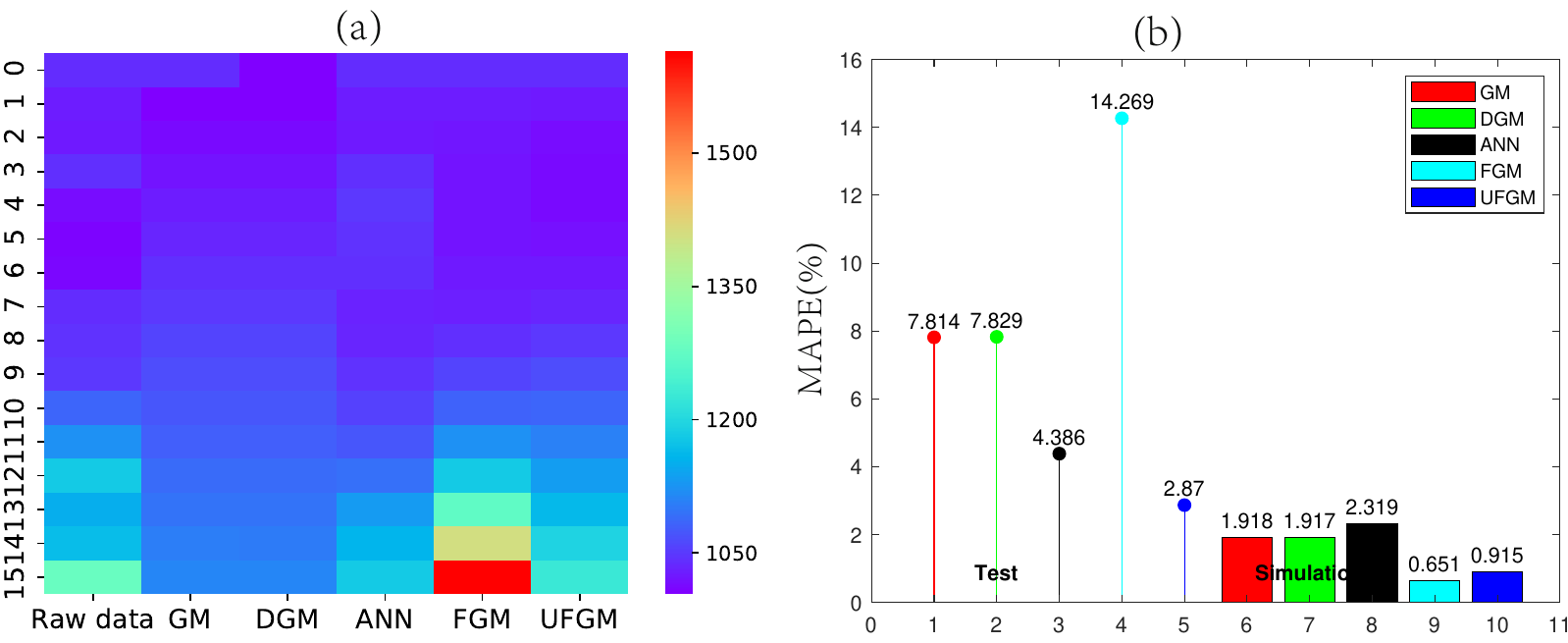}
\caption{Prediction performance of the five competitors in Case 1: (a)Simulated and predicted values by the competitors; (b)Errors by the competitors.}
\label{fig:fit1}
\end{figure}

From Table 3, it is found that the
MAPE value of UFGM(1,1) in the in-sample period is 0.91541\%, while that of competitive models, namely GM(1,1), DGM(1,1),
ANN and FGM(1,1), are 1.9183\%, 1.9171\%, 2.319511\%, and 0.65119\%, respectively, indicating the FGM(1,1) and UFGM(1,1) models have the higher accuracy in this stage because their MAPE values are closer to zeros. As for the out-of-sample period, the MAPE values of the current model and benchmarks are  2.8704\%, 7.8142\%, 7.8292\%, 4.386027\% and 14.2690\%, respectively. It is obvious that the MAPE value of the UFGM(1,1) model is far less than others, this means the UFGM(1,1) model achieves the precise prediction in this experiment. The same finding can be also obtained from Figure \ref{fig:fit1}.

\begin{figure}
\includegraphics[scale=0.52]{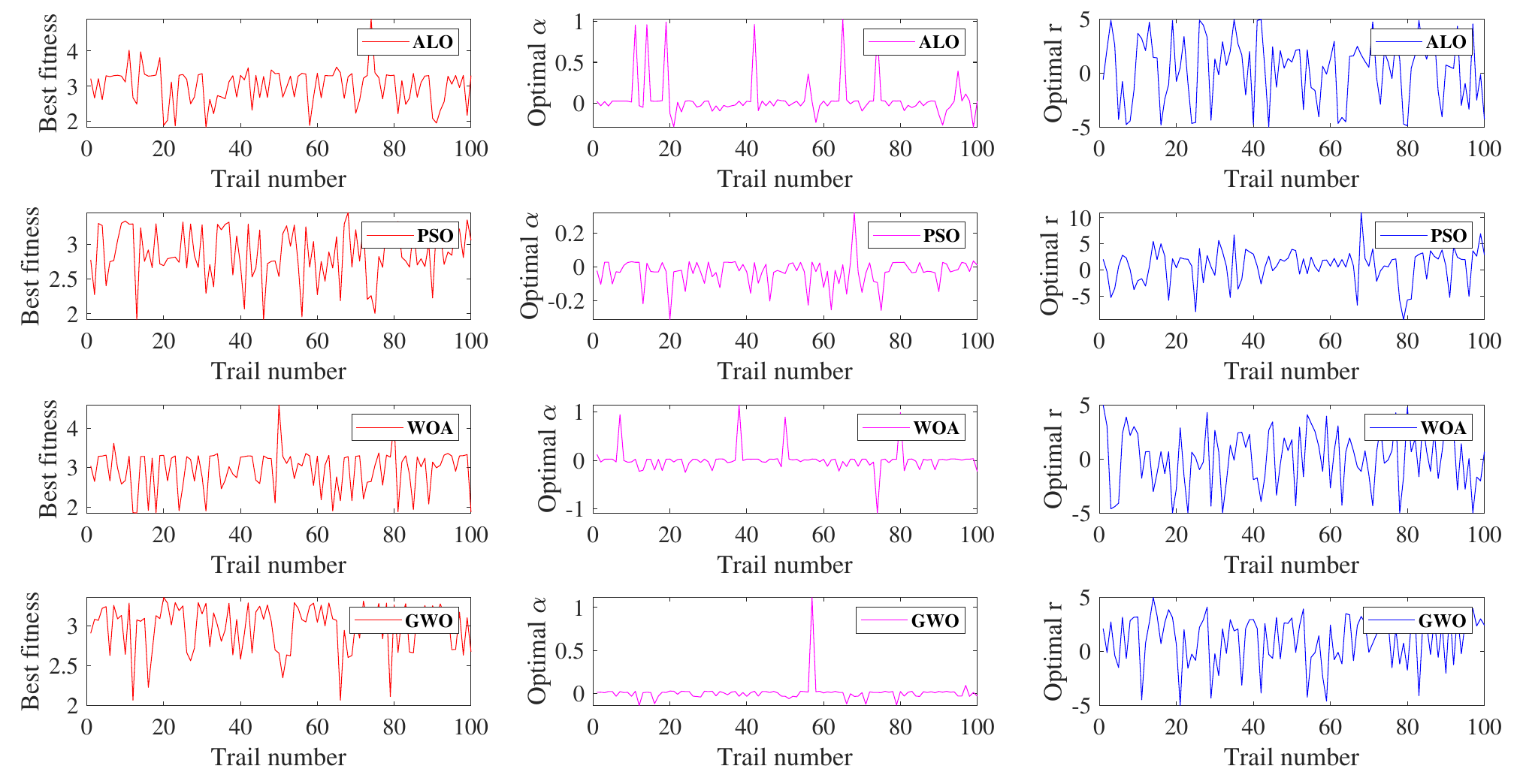}
\caption{$r$, $\alpha$ and the fitness value of the UFGM(1,1) model in each trial by using the four algorithms in Case 2.}
\label{fig:search1}
\end{figure}

\begin{figure}
\centering
\includegraphics[scale=0.6]{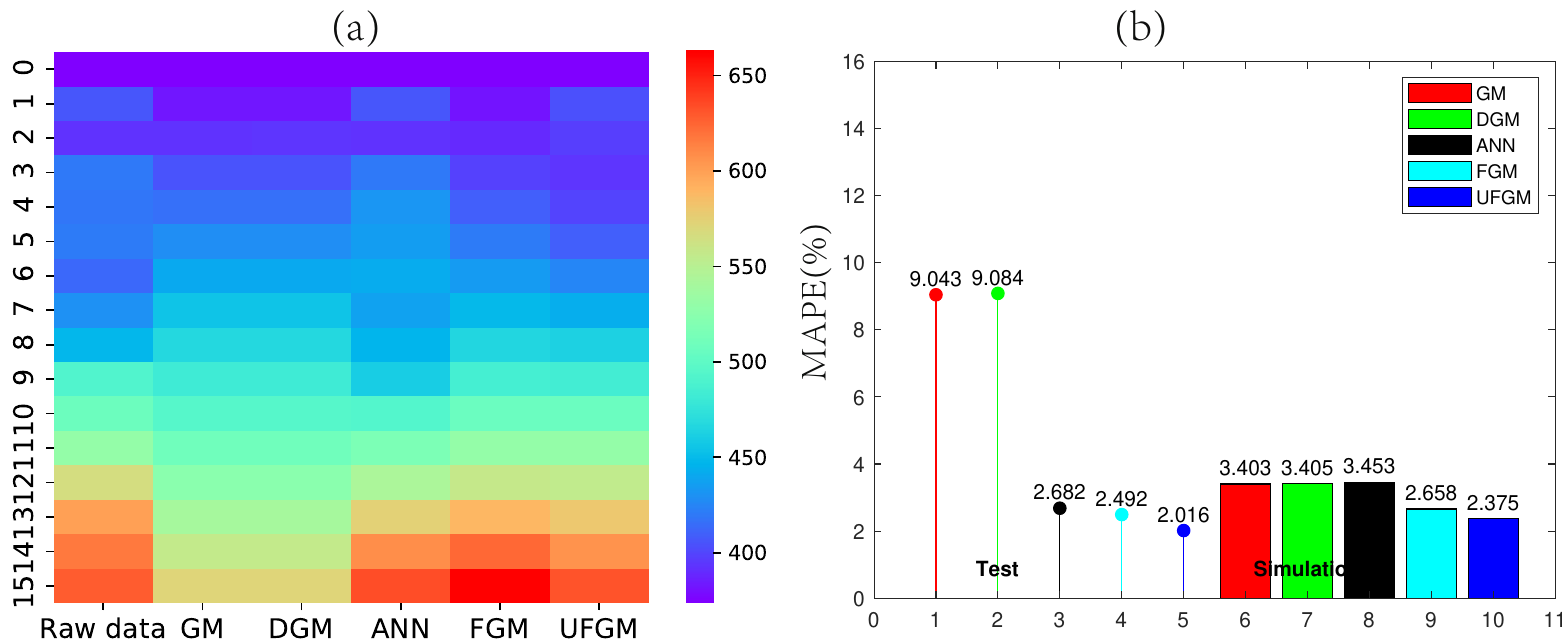}
\caption{Prediction performance of the five competitors in Case 2: (a)Simulated and predicted values by the competitors; (b)Errors by the competitors.}
\label{fig:fit2}
\end{figure}

\textbf{Case 2.} (Forecasting Chongqing's water supply production capacity.) Similar to case one, the raw data is divided into two groups, one, training set, from 2004 to 2015 are used for calibrating structural parameters; another, testing set, from 2016 to 2019 are used to evaluate the prediction performance.
First and foremost, according to the performance of the four algorithms on the determination of emerging coefficients shown in Figure \ref{fig:search1}, it is found that the fitness value by PSO is smallest, this means the emerging coefficients of UFGM(1,1) model will be determined by this algorithm. Subsequently,
the prediction results by these competitive models are obtained by reference with the current study and references therein, as shown in Table \ref{table:c2}. Meanwhile, these values in Table \ref{table:c2} are vividly displayed in Figure \ref{fig:fit2}.

\begin{table}[htbp]
  \centering
  \caption{Simulated and predicted results by the five competitive models in Case 2}\scriptsize
    \begin{tabular}{ccccccccccc}
    \hline
    Raw data & \multicolumn{1}{c}{GM(1,1)} & Error(\%) & \multicolumn{1}{c}{DGM(1,1)} & Error(\%) & \multicolumn{1}{c}{ANN} & Error(\%) & \multicolumn{1}{c}{FGM(1,1)} & Error(\%) & \multicolumn{1}{c}{UFGM} & Error(\%) \\
    \hline
    \multicolumn{11}{c}{In-sample}\\
    \hline
    373.65 & \multicolumn{1}{c}{373.65} & 0.000     & \multicolumn{1}{c}{373.65} &  0.000      & \multicolumn{1}{c}{373.65} &  0.000      & \multicolumn{1}{c}{373.65} &  0.000      & \multicolumn{1}{c}{373.65} &  0.000  \\
    405.61 & \multicolumn{1}{c}{381.59} & 5.921945 & \multicolumn{1}{c}{381.79} & 5.872636 & \multicolumn{1}{c}{405.61} &  0.000      & \multicolumn{1}{c}{380.87} & 6.099455 & \multicolumn{1}{c}{403.23} & 0.586771 \\
    391.7 & \multicolumn{1}{c}{392.78} & 0.275721 & \multicolumn{1}{c}{392.96} & 0.321675 & \multicolumn{1}{c}{391.7} &  0.000      & \multicolumn{1}{c}{389.02} & 0.684197 & \multicolumn{1}{c}{396.63} & 1.258616 \\
    419.2 & \multicolumn{1}{c}{404.31} & 3.552004 & \multicolumn{1}{c}{404.46} & 3.516221 & \multicolumn{1}{c}{419.2} &  0.000      & \multicolumn{1}{c}{398.22} & 5.004771 & \multicolumn{1}{c}{393.74} & 6.073473 \\
    418.16 & \multicolumn{1}{c}{416.17} & 0.475894 & \multicolumn{1}{c}{416.29} & 0.447197 & \multicolumn{1}{c}{431.3829} & 3.162163 & \multicolumn{1}{c}{408.6} & 2.286206 & \multicolumn{1}{c}{398.61} & 4.675244 \\
    420.35 & \multicolumn{1}{c}{428.38} & 1.910313 & \multicolumn{1}{c}{428.47} & 1.931724 & \multicolumn{1}{c}{435.0692} & 3.501649 & \multicolumn{1}{c}{420.32} & 0.007137 & \multicolumn{1}{c}{409.54} & 2.571666 \\
    412.3 & \multicolumn{1}{c}{440.94} & 6.946398 & \multicolumn{1}{c}{441} & 6.960951 & \multicolumn{1}{c}{441.3492} & 7.045642 & \multicolumn{1}{c}{433.55} & 5.154014 & \multicolumn{1}{c}{424.63} & 2.990541 \\
    429.27 & \multicolumn{1}{c}{453.88} & 5.732989 & \multicolumn{1}{c}{453.9} & 5.737648 & \multicolumn{1}{c}{437.3996} & 1.89382 & \multicolumn{1}{c}{448.47} & 4.472709 & \multicolumn{1}{c}{442.55} & 3.093624 \\
    447.83 & \multicolumn{1}{c}{467.2} & 4.325302 & \multicolumn{1}{c}{467.18} & 4.320836 & \multicolumn{1}{c}{446.0885} & 0.388882 & \multicolumn{1}{c}{465.31} & 3.903267 & \multicolumn{1}{c}{462.45} & 3.264632 \\
    491.22 & \multicolumn{1}{c}{480.9} & 2.100892 & \multicolumn{1}{c}{480.85} & 2.11107 & \multicolumn{1}{c}{459.6022} & 6.436593 & \multicolumn{1}{c}{484.32} & 1.404666 & \multicolumn{1}{c}{483.83} & 1.504418 \\
    506.89 & \multicolumn{1}{c}{495.01} & 2.343704 & \multicolumn{1}{c}{494.92} & 2.361459 & \multicolumn{1}{c}{493.1399} & 2.712648 & \multicolumn{1}{c}{505.78} & 0.218982 & \multicolumn{1}{c}{506.36} & 0.104559 \\
    529.92 & \multicolumn{1}{c}{509.53} & 3.847751 & \multicolumn{1}{c}{509.4} & 3.872283 & \multicolumn{1}{c}{516.7286} & 2.489319 & \multicolumn{1}{c}{529.99} & 0.01321 & \multicolumn{1}{c}{529.85} & 0.01321 \\
    MAPE  &       & 3.4031 &       & 3.4052 &       & 3.45384 &       & 2.6586 &       & \textbf{2.3759} \\
        \hline
    \multicolumn{11}{c}{Out-of-sample}\\
    \hline
    566.12 & \multicolumn{1}{c}{524.48} & 7.355331 & \multicolumn{1}{c}{524.3} & 7.387126 & \multicolumn{1}{c}{543.2325} & 4.042864 & \multicolumn{1}{c}{557.32} & 1.554441 & \multicolumn{1}{c}{554.15} & 2.114393 \\
    599.87 & \multicolumn{1}{c}{539.87} & 10.00217 & \multicolumn{1}{c}{539.64} & 10.04051 & \multicolumn{1}{c}{574.3076} & 4.26133 & \multicolumn{1}{c}{588.16} & 1.95209 & \multicolumn{1}{c}{579.16} & 3.452415 \\
    616.99 & \multicolumn{1}{c}{555.7} & 9.93371 & \multicolumn{1}{c}{555.43} & 9.977471 & \multicolumn{1}{c}{607.9011} & 1.473103 & \multicolumn{1}{c}{622.96} & 0.967601 & \multicolumn{1}{c}{604.82} & 1.972479 \\
    627.76 & \multicolumn{1}{c}{572.01} & 8.880782 & \multicolumn{1}{c}{571.68} & 8.93335 & \multicolumn{1}{c}{633.7528} & 0.954632 & \multicolumn{1}{c}{662.24} & 5.492545 & \multicolumn{1}{c}{631.08} & 0.528865 \\
    MAPE  &       & 9.043 &       & 9.0848 &       & 2.682982 &       & 2.4922 &       & \textbf{2.0169} \\
    \hline
    \end{tabular}%
  \label{table:c2}%
\end{table}%

Firstly, we observe from Table \ref{table:c2} and Figure \ref{fig:fit2} that for the in-sample period all of the competitors have a relatively satisfactory performance because their outputs of simulating and predicting Chongqing's water supply production capacity are close to actual ones. However, in the out-of-sample period, the predicted results by the GM(1,1) and DGM(1,1) models  deviate from the actual ones although it fall within the acceptable range. Among these models, the UFGM(1,1) model have a higher accuracy because the corresponding prediction values are closer to the actual amount of Chongqing's water supply production capacity.

Next, the error-value indicator, i.e. MAPE, is applied to further analyze the prediction performance. In the simulation stage, the MAPE values of the UFGM(1,1) and other competitive models are  2.3759\%, 3.4031\%, 3.4052\%, 3.45384\% and 2.6586\%, respectively;  those of competitors are 2.01691\%, 9.0431\%, 9.08481\%, 2.6829821\% and 2.49221\%, respectively, for the prediction stage. It is concluded that all of models are suitable for forecasting Chongqing's water supply production capacity in accordance with Lewis standard \cite{50}. To be specific, the GM(1,1) and DGM(1,1) models have a relatively poor prediction accuracy due to their higher MAPEs, the UFGM(1,1) has a superiority over others because its MAPEs are smaller either in the in-sample or out-of-sample period. It is worth mentioned that the FGM(1,1) model has a second best performance due to its second smallest MAPE. In consequence, the UFGM(1,1) will be treated as an optimal tool for forecasting Chongqing's water supply production capacity in this experiment, and the FGM(1,1) model should be considered as a sub-optimal technique in dealing with Chongqing's water supply production capacity.

\section{Conclusion and summary}
\label{sec:6}
In this study, an unified framework for fractional grey models is presented by combining general fractional derivative with memory effects and grey modeling theory. This unified framework (denoted as UFGM(1,1) for short) not just covers commonly-used fractional grey models already in place, but also can derived other new ones. The main conclusions are summarized as follows.
\begin{itemize}
\item Based on the GR and GRL fractional derivatives, the framework for published fractional grey models is developed. By taking different kernel functions and normalization ones, some other new fractional grey models can be deduced, a new fractional grey model, for example, with $k(x,r) = e^{-\frac{r}{1 - r}x}$ and $N(r)=1$, is expressed as $
_0^{\exp }D_t^r{x^{(\alpha )}}(t) + a{x^{(\alpha )}}(t) = bt + c.
$
\item The commonly-used intelligent algorithms, namely GWO, PSO, WOA and ALO, are employed to search for the optimal emerging coefficients, i.e., $r$, and $\alpha$, for the UFGM(1,1) model in the current study.
    \item Two published cases are utilized to verify the superiority of the UFGM(1,1) model and explore the effects of fractional accumulation order and initial value on the prediction performance.
        \item The UFGM(1,1) model is also applied to dealing with Henan's and Chongqing's water supply production capacity so as to explain the validity of the proposed method. Numerical results show that the UFGM(1,1) model outperforms a range of benchmarks.
\end{itemize}

In addition, the unified framework is not without limitations, for example, this framework is focusing on continuous-fractional grey models, potentially neglecting the discrete-fractional grey model (e.g. conformable fractional grey model \cite{24}), this is a primary research direction of our following study. As stated in \cite{51}, the data-preprocessing techniques (e.g., data grouping \cite{52}, rolling mechanism \cite{53}, etc.) should be also worth exploring in the coming work.

\noindent\textbf{Acknowledgements}

This work was supported by the Postgraduate Research \& Practice Innovation
Program of Jiangsu Province (KYCX20\_1144) and the Fundamental Research Funds
for the Central Universities of China (Grant no. 2019YBZZ062).

\end{document}